[1996/12/01]
\documentclass{article}
\usepackage{amsmath,amsthm}
\newtheorem{thm}{Theorem}[section]

\theoremstyle{definition}

\newcommand{\floor}[1]{\left\lfloor#1\right\rfloor}

\theoremstyle{remark}

\title{On Sinha's note on perfect numbers
}

\date{}
\author{Tomohiro Yamada}
\begin{document}
\maketitle

\begin{abstract}
We shall show that there is no odd perfect number of the form $2^n+1$ or $n^n+1$.
\end{abstract}

\section{Introduction}\label{intro}
A positive integer $N$ is called perfect if $\sigma(N)=2N$,
where $\sigma(N)$ denotes the sum of divisors of $N$.
As is well known, an even integer $N$ is perfect if and only if
$N=2^{k-1}(2^k-1)$ with $2^k-1$ prime.
In contrast, one of the oldest unsolved problems is whether there exists
an odd perfect number or not.  Moreover, it is also unknown whether
there exists an odd $m$-perfect number for an integer $m\geq 2$, i.e.,
an integer $N$ with $\sigma(N)=mN$ or not.

Sinha \cite{Sin} showed that $28$ is the only even perfect number of the form
$x^n+y^n$ with $\gcd(x, y)=1$ and $n\geq 2$ and also the only even perfect number
of the form $a^n+1$ with $n\geq 2$.
On the other hand,
it is not even proved or disproved that there exists no odd perfect number of the form
$x^2+1$ with $x$ an integer.
Klurman \cite{Klu} proved that if $P(x)$ is a polynomial of degree $\geq 3$ without repeated factors,
then there exist only finitely many odd perfect numbers of the form $P(x)$ with $x$ an integer.
Luca \cite{Luc} (cited in Theorem 9.8 of \cite{KLS}) showed that no Fermat number can be perfect.

In this article, we would like to prove that there exists no odd perfect number
of the form $2^n+1$ or $n^n+1$.

Indeed, we prove a more general result.

\begin{thm}\label{thm1}
Let $m$ and $U$ be nonnegative integers.
We put \\ $s_0=\floor{2^U\log a/(U+1)\log 2}$ and
$t_0=2s_0+1$ if $U=0$ and $a+1$ is square and $t_0=2s_0$ otherwise.
Let $c=1.093\cdots=(\log 2)/2+(\log 3)/3-(\log^2 3)/2$
and $C=C(U)$ be the constant defined by
\[C=\sum_{2^{U+1} (2m+1)<16, } \frac{1-\log\log(2^{U+1} m)}{2^{U+1} m}.\]

If $a^n+1$ is an odd $(4m+2)$-perfect number and $n=2^U$, then
\begin{equation}\label{eq1}
\log a>\frac{((4m+2)/e^C)^{2^{U+1}}}{2^U}.
\end{equation}
If $a^n+1$ is an odd $(4m+2)$-perfect number and $n=2^U v$ with $v>1$ odd, then
\begin{equation}\label{eq2}
\begin{split}
& \log (4m+2)-C \\
< & \frac{\exp\left(\frac{1+\log t_0}{2^{U+1}}\right)}{2^{U+1}}
\left(\log(2^U \log a)+(U+1)(1+\log t_0)\log 2 +\frac{\log^2 t_0}{2}+c\right).
\end{split}
\end{equation}
Moreover, no integer of the form $2^n+1$ can be $(4m+2)$-perfect.
\end{thm}

For example, if $a^{128s}+1$ is odd $(4m+2)$-perfect, then $a\geq 10$
and, if $a^{256s}+1$ is odd $(4m+2)$-perfect, then $a\geq 18$.
Furthermore, if $a^{16}+1$ is odd $(4m+2)$-perfect, then $a>\exp\exp 19.4$
and, if $a^{32}+1$ is odd $(4m+2)$-perfect, then $a>\exp\exp 40.8$.
We note that $C(0)=0.9807\cdots$, $C(1)=0.1758\cdots$, $C(2)=0.03348\cdots$
and $C(U)=0$ for $U\geq 3$.

We shall prove that an odd perfect number of the form $n^n+1$
must be of the form $2^m+1$
and deduce the following result from the above result.

\begin{thm}\label{thm2}
$28$ is the only $(4m+2)$-perfect number of the form $n^n+1$ with $m, n\geq 0$ an integer.
\end{thm}

Thus, we conclude that $28$ is the only perfect number of the form $n^n+1$.

\section{Proof of Theorem \ref{thm1}}

Assume that $a^n+1$ is an odd $(4m+2)$-perfect number.
By Euler's result, we must have $a^n+1=px^2$ for a prime $p$ and an integer $x$.

Write $n=2^U p_1^{e_1} p_2^{e_2} \ldots p_r^{e_r}$ with $p_1>p_2>\cdots >p_r$ odd primes
and let $P_i=p_i^{e_i}$ for $i=1, 2, \ldots, r$ and $s=\omega(a^{2^U}+1)$.
We put $o_p(x)$ to be the multiplicative order of $x$ modulo $p$.

We can factor $a^n+1=M_0 M_1 \cdots M_r$,
where $M_0=a^{2^U}+1$ and
\[M_i=\frac{a^{2^U P_1 P_2 \cdots P_i}+1}{a^{2^U P_1 P_2 \cdots P_{i-1}}+1}\]
for $i=1, 2, \ldots, r$.
Moreover, let \[L_i=M_0 M_1 \ldots M_i=a^{2^U P_1 P_2 \cdots P_i}+1\]
and $M_i=E_i Y_i^2, L_i=D_i X_i^2$ with $D_i$ and $E_i$ squarefree.
Clearly, we have $a^n+1=L_r=px^2$ and therefore $D_r=p$.

We begin by showing that $p_i\equiv 1\pmod{2^{U+1}}$ for every $i$.
If $\gcd((a^n+1)/(a^{n/P_i}+1), a^{n/P_i}+1)=1$, then
\begin{equation}\label{eqa}
a^{n/P_i}+1=X^2, \frac{a^n+1}{a^{n/P_i}+1}=pY^2
\end{equation}
or
\begin{equation}\label{eqb}
a^{n/P_i}+1=pX^2, \frac{a^n+1}{a^{n/P_i}+1}=Y^2
\end{equation}
for some integers $X$ and $Y$.
If $U=0$, then we clearly have $p_i\equiv 1\pmod{2^{U+1}}$.
If $U>0$, then $n/p_i^{e_i}$ is even and \eqref{eqa} is clearly impossible.
The impossibility of \eqref{eqb} follows from Ljunggren's result \cite{Lju} that
$(a^f+1)/(a+1)$ with $a\geq 2, f\geq 3$ cannot be square.

Hence, we must have $\gcd((a^n+1)/(a^{n/P_i}+1), a^{n/P_i}+1)>1$.
Observing that
\[\frac{a^n+1}{a^{n/P_i}+1}=\sum_{j=0}^{P_i-1}(-1)^j a^{j(n/P_i)}\equiv P_i\pmod{a^{n/P_i}+1},\]
$p_i$ must divide $a^{n/P_i}+1$.
Thus, proceeding as in the proof of Theorem 4.12 of \cite{KLS},
we see that $2^{U+1}$ divides $o_{p_i}(a)$ and $o_{p_i}(a)$ divides $2n/P_i$.
In particular, $p_i\equiv 1\pmod{2^{U+1}}$ for every $i$.

Nextly, we show that for each $i=1, 2, \ldots, r$, we have either
(i) $\gcd(L_{i-1}, M_i)=1$ and $\omega(D_{i-1})<\omega(D_i)$ or
(ii) $p_i$ is the only prime dividing $\gcd(L_{i-1}, M_i)$ and $p_i$ divides $a^{2^U}+1$.

If $\gcd(L_{i-1}, M_i)=1$, then we must have $D_i=D_{i-1} E_{i-1}$ and $X_i=X_{i-1} Y_{i-1}$.
It follows from Ljunggren's result mentioned above that $E_{i-1}\neq 1$.
Since $D_i$ is squarefree, we have $\omega(D_{i-1})<\omega(D_i)$.

Assume that $\gcd(L_{i-1}, M_i)>1$.
Since
\[M_i=\sum_{j=0}^{P_i-1} (-1)^j 2^{2^U P_1 P_2 \ldots P_{i-1}j}\equiv P_i\pmod{L_{i-1}},\]
we see that $p_i$ is the only prime dividing both $L_{i-1}$ and $M_i$.

Now $p_i$ must divide $L_{i-1}$ and therefore,
proceeding as above, we see that $2^{U+1}$ divides $o_{p_i}(a)$ and $o_{p_i}(a)$ divides $2^{U+1} P_1 P_2 \cdots P_{i-1}$.
Hence, $o_{p_i}(a)=2^{U+1} d$ and therefore $p_i\equiv 1\pmod{2^{U+1} d}$ for some $d$ dividing $P_1 P_2 \cdots P_{i-1}$.
But, since $p_1>\cdots >p_{i-1}>p_i$, we must have $o_{p_i}(a)=2^{U+1}$ and therefore
$p_i$ must divide $a^{2^U}+1$.

It is clear that (ii) occurs at most $s$ times.
Moreover, we observe that in the case (ii), $p_i$ is the only possible prime
which divides $D_{i-1}$ but not $D_i$.
Hence, we must have $\omega(D_{i-1})\leq \omega(D_i)+1$ for each $i$.
Now we see that (i) also occurs at most $s$ times.

We can easily see that $\omega(D_0)=0$ if and only if $U=0$ and $a+1$ is a square.
Thus we conclude that $r\leq 2s+1$ if $D_0=a+1$ with $U=0$ is square and $r\leq 2s$ otherwise.

If a prime $p$ divides $a^{2^U d}+1$ but $a^{2^U e}+1$ for any $e<d$,
then the multiplicative order of $2\pmod{p}$ is equal to $2^{U+1} d$ and therefore
$p=2^{U+1}kd+1$ for some integer $k$.
Moreover, the number of such primes is at most $k_0(d)=\floor{2^U d\log a/\log (2^{U+1} d)}$
and therefore $s\leq s_0$.

Hence, for each $d$,
\begin{equation}
\begin{split}
\prod_{o_p(a)=2^{U+1} d}\frac{p}{p-1}< & \exp \sum_{o_p(a)=2^{U+1} d}\frac{1}{p-1}\leq \sum_{k=1}^{k_0(d)} \frac{1}{2^{U+1}kd} \\
\leq & \exp\frac{1+\log(2^U d\log a/\log(2^{U+1} d))}{2^{U+1} d},
\end{split}
\end{equation}
so that
\begin{equation}
\frac{\sigma(a^n+1)}{a^n+1}=\prod_{\substack{o_p(a)=2^{U+1} d,\\ d\mid P_1 P_2 \ldots P_r}}\frac{p}{p-1}
<\exp\left(C+\sum_{d\mid P_1 P_2 \ldots P_r}\frac{\log (2^U d\log a)}{2^{U+1} d}\right).
\end{equation}

If $r=0$, then we immediately see that
\begin{equation}
\sum_{d\mid P_1 P_2 \ldots P_r}\frac{\log(2^U d\log a)}{2^{U+1} d}=\frac{U\log 2+\log\log a}{2^{U+1}}.
\end{equation}

If $r>0$, then, observing that
\begin{equation}
\sum_{i=0}^\infty \frac{i}{q^i}=\sum_{j=0}^\infty \sum_{i=j+1}^\infty \frac{1}{q^i}=\sum_{j=0}^\infty\frac{1}{q^j(q-1)}=\frac{q}{(q-1)^2},
\end{equation}
we have
\begin{equation}\label{eqc}
\begin{split}
& \sum_{d\mid P_1 P_2 \ldots P_r}\frac{\log(2^U d\log a)}{2^{U+1} d} \\
< & \sum_{f_1, f_2, \ldots, f_r\geq 0} \frac{\log (2^U\log a)+f_1\log p_1+f_2\log p_2+\cdots +f_r \log p_r}{2^{U+1} p_1^{f_1} p_2^{f_2} \cdots p_r^{f_r}} \\
= & \prod_{i=1}^t\frac{p_i}{p_i-1}\left(\frac{\log (2^U\log a)}{2^{U+1}}+\sum_{k=1}^t \frac{\log p_k}{2^{U+1} (p_k-1)}\right) \\
= & \left(\frac{1}{2^{U+1}}\prod_{i=1}^r\frac{p_i}{p_i-1}\right)\left(\frac{\log (2^U\log a)}{2^{U+1}}+\sum_{k=1}^r \frac{\log p_k}{p_k-1}\right).
\end{split}
\end{equation}

Since each $p_i\equiv 1\pmod {2^{U+1}}$,
we have
\begin{equation}
\prod_{i=1}^r \frac{p_i}{p_i-1}<\prod_{k=1}^r \frac{2^{U+1} k+1}{2^{U+1} k}<\exp \frac{1+\log r}{2^{U+1}}
\end{equation}
and observing that $\sum_{k=1}^t \log k/k\leq (\log t)^2/2+c$ for $t\geq 1$,
\begin{equation}
\begin{split}
\sum_{k=1}^r \frac{\log p_k}{p_k-1}
< & \sum_{k=1}^r \frac{\log k+(U+1)\log 2}{2^{U+1} k} \\
< & \frac{1}{2^{U+1}}\left((U+1)(1+\log r)\log 2 +\frac{\log^2 r}{2}+c\right).
\end{split}
\end{equation}
Thus, we obtain
\begin{equation}
\begin{split}
& \sum_{d\mid P_1 P_2 \ldots P_r}\frac{\log(2^U d\log a)}{2^{U+1} d} \\
< & \frac{\exp\left(\frac{1+\log r}{2^{U+1}}\right)}{2^{U+1}}\left(\log(2^U \log a)+(U+1)(1+\log r)\log 2 +\frac{\log^2 r}{2}+c\right).
\end{split}
\end{equation}

We see that $r\leq t_0$, where we recall that $s\leq s_0=\floor{2^U\log a/(U+1)\log 2}$.
Hence, we conclude that
\begin{equation}
\log(4m+2)=\log\frac{\sigma(a^n+1)}{a^n+1}<C+\frac{U\log 2+\log\log a}{2^{U+1}}
\end{equation}
if $r=0$ and
\begin{equation}
\begin{split}
& \log (4m+2)-C \\
< & \frac{\exp\left(\frac{1+\log t_0}{2^{U+1}}\right)}{2^{U+1}}
\left(\log(2^U \log a)+(U+1)(1+\log t_0)\log 2 +\frac{\log^2 t_0}{2}+c\right)
\end{split}
\end{equation}
otherwise.
Thus \eqref{eq1} and \eqref{eq2} follows.

Now we consider the case $a=2$.
If $U\geq 4$, then the right-hand side of \eqref{eq1} and \eqref{eq2} is $<0.53<\log 2$
and therefore $a^n+1$ cannot be $(4m+2)$-perfect.

If $U\leq 3$, then $2^{2^U}+1$ is prime and therefore $s=1$.
Clearly, for $n=2^U$ with $U\leq 3$, $2^n+1=2^{2^U}+1$ is not $(4m+2)$-perfect.
Hence, we must have $r\leq 2$ and $n=2^U p_1^{e_1}$ or $2^U p_1^{e_1} p_2^{e_2}$.

If $n=2^U p_1^{e_1}$, then, iterating the argument given before,
we must have $p_1=2^{2^U}+1$.
Thus, $n=3^{e_1}$, $2\times 5^{e_1}$, $2^2\times 17^{e_1}$ or $2^3\times 257^{e_1}$.

However, for $n=3^{e_1}$ with $e_1\geq 3$, we see that both primes $19$ and $87211$
divide $2^n+1$ exactly once
since $19$ and $87211$ divide $2^{27}+1$ exactly once and the only prime dividing
both $(2^n+1)/(2^{27}+1)$ and $2^{27}+1$ is $3$.
This implies that $2^n+1$ cannot be of the form $px^2$ and therefore $2^n+1$ cannot be $(4m+2)$-perfect
if $n=3^{e_1}$ with $e_1\geq 3$.
Similarly, $41$ and $101$ divide $2^n+1$ exactly once if $n=2\times 5^{e_1}$ and $e_1\geq 2$.
Clearly, none of $2^3+1, 2^9+1, 2^{10}+1$ is $(4m+2)$-perfect.
Thus $2^n+1$ cannot be $(4m+2)$-perfect if $n=3^{e_1}$ or $2\times 5^{e_1}$.
Similarly, $2^n+1$ cannot be $(4m+2)$-perfect if $n=2^2 \times 17^{e_1}$ or $2^3 \times 257^{e_1}$.

If $n=2^U p_1^{e_1}p_2^{e_2}$, then, iterating the argument given before,
$p_1>p_2=2^{2^U}+1$.

If $U=1$ and $n=10 p_1^{e_1}$, then we must have
\[2^{10}+1=5^2\times 41, \frac{2^n+1}{2^{10}+1}=41py^2\]
since $(2^n+1)/(2^{10}+1)$ cannot be square by Ljunggren's result.
Thus, we must have $p_1=41$.
However, this implies that $2^n+1$ must be divisible by $821$ and $10169$ exactly once,
which contradicts to the fact that $2^n+1=px^2$.
If $U=1$ and $n=2\times 5^{e_2} p_1^{e_1}$ with $e_2\geq 2$, then,
since three primes $41, 101, 8101$ divide $2^{50}+1$ exactly once,
at least two of these primes divide $2^n+1$.
Thus $2^n+1$ cannot be $(4m+2)$-perfect if $n=2p_1^{e_1}p_2^{e_2}$.
Similarly, $2^n+1$ cannot be $(4m+2)$-perfect for $n=2^U p_1^{e_1}p_2^{e_2}$ with $U=2, 3$.

Now we assume that $n=3^{e_2} p_1^{e_1}$.

If $n=3^{e_2} p_1^{e_1}$ with $e_2\geq 4$, then,
at least two of three primes $19, 163, 87211$ divide $2^n+1$ exactly once and therefore
$2^n+1$ cannot be $(4m+2)$-perfect for such $n$.
If $n=27p_1^{e_1}$, then we must have $p_1=19$ or $87211$.
We cannot have $p_1=19$ since $571$ and $87211$ divide $2^n+1$ exactly once for $n=27\times 19^{e_1}$.
Assume that $p_1=87211$.
We observe that, for $d=3^{f_2} 87211^{f_1}$ with $f_1>0$, we have
\begin{equation}
\prod_{o_p(a)=2d}\frac{p}{p-1}<\exp \frac{1+\log(d\log 2/\log(2d))}{2d}<\exp \frac{\log d}{2d}
\end{equation}
and, proceeding as in \eqref{eqc},
\begin{equation}\label{eqd}
\sum_{\substack{d=3^{f_2} 87211^{f_1},\\ f_1>0, f_2\geq 0}} \frac{\log d}{2d}<\frac{87211}{116280}\left(\frac{\log 3}{174422}+\frac{\log 87211}{87210}\right)<\frac{1}{9000}.
\end{equation}
Thus, $\sigma(2^n+1)/(2^n+1)<e^{1/9000}\sigma(2^{27}+1)/(2^{27}+1)<2$
and therefore $2^n+1$ cannot be $(4m+2)$-perfect.

If $n=9p_1^{e_1}$, then we must have $p_1=19$
and therefore two primes $571$ and $174763$ divide $2^n+1$ exactly once,
which is a contradiction.

Finally, assume that $n=3p_1^{e_1}$.
If $p_1\geq 11$, then, like \eqref{eqd},
\begin{equation}
\sum_{\substack{d=3^{f_2} p_1^{f_1},\\ f_1>0, f_2\geq 0}} \frac{\log d}{2d}<\frac{3p_1}{2(p_1-1)}\left(\frac{\log 3}{2p_1}+\frac{\log p_1}{p_1-1}\right)<0.24
\end{equation}
and $\sigma(2^n+1)/(2^n+1)<(13/9)e^{0.24}<2$, which is a contradiction.

The only remaining case is $n=3p_1^{e_1}$ with $p_1=5$ or $7$.
We observe that $2^{15}+1=3^2\times 11\times 331$ and $2^{21}+1=3^2\times 43\times 5419$.
Thus $2^n+1$ must be divisible by at least two distinct primes exactly once,
which is a contradiction again.
Now we conclude that $2^n+1$ can never be $(4m+2)$-perfect.

\section{Proof of Theorem \ref{thm2}}

Sinha's result clearly implies that $28$ is the only even perfect number of the form $n^n+1$.
Thus, we may assume that $n^n+1$ is an odd $(4m+2)$-perfect number.
Clearly $n$ must be even and we can write $n=2^u s$ with $u>0$ and $s$ odd.

As before, we must have
$n^n+1=px^2$ for some prime $p$ and integer $x$.

Assume that $s>1$.
Then we must have
\begin{equation}
n^n+1=(n^{2^u}+1)\times \frac{n^{2^u s}+1}{n^{2^u}+1}=N_1 N_2,
\end{equation}
say.

If $N_1$ and $N_2$
have a common prime factor $p$, then $p$ divides $d_2$ and therefore $p$ divides $2^u s=n$.
This is impossible since $\gcd(n^n+1, n)=1$.
Thus, we see that $\gcd(N_1, N_2)=1$ and therefore
$N_1=X^2, N_2=pY^2$ or $N_1=pX^2, N_2=Y^2$.

We can easily see that $n^{2^u}+1$ cannot be square since $u>0$
and therefore
\begin{equation}
\frac{n^{2^u s}+1}{n^{2^u}+1}=Z^2.
\end{equation}
However, this is also impossible from Ljunggren's result.

Now we must have $s=1$ and $n^n+1=2^{u 2^u}+1$,
which we have just proved not to be $(4m+2)$-perfect in Theorem \ref{thm1}.
This proves Theorem \ref{thm2}.

{}
\vskip 12pt

{\small Tomohiro Yamada}\\
{\small Center for Japanese language and culture\\Osaka University\\562-8558\\8-1-1, Aomatanihigashi, Minoo, Osaka\\Japan}\\
{\small e-mail: \protect\normalfont\ttfamily{tyamada1093@gmail.com}}
\end{document}